\documentclass[journal,twoside,web]{ieeecolor}
\usepackage{lcsys}
\usepackage{cite}
\usepackage{amsmath,amssymb,amsfonts}
\usepackage{algorithmic}
\usepackage{graphicx}
\usepackage{textcomp}
\usepackage{hhline}

\newtheorem{theorem}{Theorem}
\newtheorem{assumption}{Assumption}

\usepackage{eso-pic}
\AddToShipoutPictureBG*{%
\AtPageUpperLeft{%
\setlength\unitlength{1in}%
\hspace*{\dimexpr0.5\paperwidth\relax}
\makebox(0,-0.5)[c]{
\begin{tabular}{c c}
To appear in the IEEE Control Systems Letters.
Uploaded to ArXiv \today. \\
\end{tabular}}}}

\AddToShipoutPictureBG*{%
\AtPageUpperLeft{%
\setlength\unitlength{1in}%
\hspace*{\dimexpr0.5\paperwidth\relax}
\makebox(0,-21.4)[c]{
\footnotesize
\begin{tabular}{c c}
© 2025 IEEE. Personal use of this material is permitted. Permission from
IEEE must be obtained for all other uses, in any \\
current or future media, including reprinting/republishing
this material for advertising or promotional purposes, creating new \\
collective works, for resale or redistribution to servers or lists,
or reuse of any copyrighted component of this work in other works.
\end{tabular}}}}

\def\BibTeX{{\rm B\kern-.05em{\sc i\kern-.025em b}\kern-.08em
    T\kern-.1667em\lower.7ex\hbox{E}\kern-.125emX}}
\begin{document}
\title{A Contingency Model Predictive Control Framework for Safe Learning}
\author{Merlijne Geurts$^*$, Tren Baltussen$^*$, \IEEEmembership{Member, IEEE}, Alexander Katriniok, \IEEEmembership{Senior Member, IEEE}, \newline Maurice Heemels, \IEEEmembership{Fellow, IEEE} \vspace{-10mm}
\thanks{Received March 17, 2025; revised May 5, 2025; accepted
May 19, 2025. Date of publication XXX XX 2025, date of
current version XXX XX 2025. This research has received funding from the Dutch Research Council (NWO) via AMADeuS, project no. 18489.}
\thanks{
$^*$Merlijne Geurts and Tren Baltussen contributed equally to this work. (Corresponding author: Tren Baltussen). 
}
\thanks{All authors are with the Control Systems Technology Section, Eindhoven University of Technology, The Netherlands {\tt{ \{m.e.geurts, t.m.j.t.baltussen, a.katriniok, m.heemels\}@tue.nl}}.}
\thanks{Digital Object Identifier 10.1109/LCSYS.2025.3575191}%
}

\maketitle

\begin{abstract}
This research introduces a multi-horizon contingency model predictive control (CMPC) framework in which classes of robust MPC (RMPC) algorithms are combined with classes of learning-based MPC (LB-MPC) algorithms to enable safe learning.
We prove that the CMPC framework inherits the robust recursive feasibility properties of the underlying RMPC scheme, thereby ensuring safety of the CMPC in the sense of constraint satisfaction.
The CMPC leverages the LB-MPC to safely learn the unmodeled dynamics to reduce conservatism and improve performance compared to standalone RMPC schemes, which are conservative in nature.
In addition, we present an implementation of the CMPC framework that combines a particular RMPC and a Gaussian Process MPC scheme.
A simulation study on automated lane merging demonstrates the advantages of our general CMPC framework.
\end{abstract}


\vspace{-6pt}
\section{Introduction}
\label{sec:introduction}
\IEEEPARstart{M}{odel} predictive control (MPC) is a well-established control method for safety-critical systems, as MPC can ensure safety in terms of constraint satisfaction while attaining high performance.
In particular, robust MPC (RMPC) can provide rigorous safety guarantees under all modeled uncertainties through robust recursive feasibility properties \cite{bib:alvarado_limon,bib:nezami}. Unfortunately, RMPC can be overly conservative and inflexible as it accounts for the worst-case disturbance at all times, which limits its performance and applicability.
Alternatively, learning-based MPC (LB-MPC) methods use data-driven and machine learning techniques, for instance, to improve the system prediction model, and consequently to increase performance \cite{bib:hewing}, \cite{bib:hewing2}. Although LB-MPC methods can account for uncertainty, providing rigorous safety guarantees for LB-MPC methods is still an open research problem \cite{bib:hewing2}.

In order to ensure constraint satisfaction, LB-MPC methods are often combined with robust control methods.
For example, \cite{bib:zanon} uses reinforcement learning to safely adapt the parametrization of an RMPC.
Further, robust predictive safety filters can be used to correct/modify control actions when a learning-based controller will jeopardize safety \cite{bib:wabersich}.
In \cite{bib:tomlin} an LB-MPC optimizes a learning-based cost function subject to robust constraints.
Various existing methods compute a separate high-performance plan and a robust backup plan, and use this backup plan when safety is compromised (see e.g., \cite{bib:wabersich, bib:brudigam, bib:li_bastani}). A potential drawback of these approaches is erratic intervention and switching between controllers, and possibly sub-optimal control performance. 
We refer to \cite{bib:hewing2} and \cite{bib:mesbah} for an overview on safe LB-MPC.

An alternative approach is Contingency MPC (CMPC) \cite{bib:alsterda}, which uses two separate MPC horizons to compute both a performance and a contingency input sequence. Instead of computing a single input sequence as, e.g., in \cite{bib:tomlin, bib:soloperto}, CMPC simultaneously computes two input sequences that are coupled by their first input.
As a result, CMPC can leverage the assertive nature of LB-MPC while ensuring (safety) constraint satisfaction via the underlying RMPC, without switching between the two plans.
Interestingly, there exist several relevant works in which CMPC is applied to various safety-critical systems, mostly in automated driving and traffic applications
\cite{bib:alsterda}, \cite{bib:schweidel, bib:dallas, bib:chen}.
However, these works either lack formal safety guarantees \cite{bib:alsterda, bib:schweidel, bib:dallas} or do not consider LB-MPC \cite{bib:chen}.
As CMPC presents a powerful methodology for safe learning-based MPC, a general conceptual framework with high-level conditions for robust recursive feasibility can be of great interest for safe learning.

To fill this gap and as the main contribution of this paper, we propose a conceptual CMPC framework which allows the integration of various classes of LB-MPC and RMPC schemes (Sec. \ref{sec:framework}), resulting in a framework of which \cite{bib:tomlin} and \cite{bib:chen} are special cases.
We show that the CMPC inherits the safety guarantees from the underlying RMPC by formally proving robust recursive feasibility of the CMPC framework (Sec. \ref{sec:safety}).
This will highlight the basic assumptions required by the underlying MPC schemes to enable safe learning.
Moreover, we demonstrate how the CMPC framework can be used by presenting a specific implementation using the RMPC in \cite{bib:geurts_robust} and a Gaussian Process MPC (GP-MPC) \cite{bib:baltussen} (Sec. V-VI). 
Finally, the CMPC is applied to a lane-merging use case with an autonomous vehicle and is compared to the baseline RMPC and GP-MPC (Sec. VII). CMPC shows improved performance over the standalone RMPC while guaranteeing safety as opposed to the standalone GP-MPC.

\section{Problem Formulation}
We consider the discrete-time non-linear system $\mathcal{S}$ subject to uncertainty, of the form
\begin{equation}\label{eq:system_S}
    x_{k+1} = f(x_k, u_k) + \underbrace{ g(x_k,u_k) +  v_k}_{=:w_k},
\end{equation}
where $x_k \in \mathbb{R}^{n_x}$ and $u_k\in \mathbb{R}^{n_u}$ denote the state and control input, respectively, of $\mathcal{S}$ at time step $k\in\mathbb{N} :=\{0,1,2,\dots\}$. 
The dynamics of $\mathcal{S}$ are composed of known, nominal dynamics denoted by
$f:\mathbb{R}^{n_x}\times\mathbb{R}^{n_u}\rightarrow\mathbb{R}^{n_x}$, of unknown, unmodeled dynamics $g:\mathbb{R}^{n_x}\times\mathbb{R}^{n_u}\rightarrow\mathbb{R}^{n_x}$, and process noise $v_k \in \mathbb{V} \subset \mathbb{R}^{n_x}$. 
The system is subject to state constraints $x_k\in \mathbb{X}$ and input constraints $u_k \in \mathbb{U}$, $k\in\mathbb{N}$. The disturbance $w_k= g(x_k,u_k)  + v_k$ is assumed to satisfy $w_k \in \mathbb{W} \subset \mathbb{R}^{n_x}$, for all $k \in \mathbb{N}$, where $\mathbb{W}$ is a bounded set and $0\in\mathbb{W}$.
In this paper, we aim to address the problem of controlling system $\mathcal{S}$ and learning the unmodeled dynamics $g$ to improve the performance while providing theoretical safety guarantees in the sense of constraint satisfaction.

\vspace{-2pt}
\section{CMPC Framework for Safe Learning}\label{sec:framework}
To address the problem formulated above, we propose a general multi-horizon CMPC framework that uses two prediction models, based on the following optimal control problem (OCP) for a given state $x_k$ at time $k \in \mathbb{N}$:
\begin{subequations}
\vspace{-1pt}
\label{eq:MPC_problem}
\begin{align}
  \min_{\bar{x}_{k \mid k}, \bar{U}_k,\hat{U}_k}&~ J(x_k, \bar{x}_{k \mid k}, \bar{U}_k,\hat{U}_k)\label{subeq:cost}\\
  \hspace{-1.8mm}\text{s.t. }  \bar{x}_{k+j+1 \mid k} & = f(\bar{x}_{k+j \mid k},\bar{u}_{k+j \mid k}),\  j=0,1,\dots,N-1, \label{subeq:fr} \\
  \hat{x}_{k+j+1 \mid k} & = f(\hat{x}_{k+j \mid k},\hat{u}_{k+j \mid k}) + d(\hat{x}_{k+j \mid k},\hat{u}_{k+j \mid k}, \mathcal{D}_k) \notag \\  & \quad \quad j=0,1,\dots,N-1,\label{subeq:fL}  \\
  \bar{x}_{k \mid k} & = \hat{x}_{k \mid k} \in x_{k} \oplus X
  , \label{subeq:same_initial_x}\\
  \bar{u}_{k \mid k} & = \hat{u}_{k \mid k}, \label{subeq:same_initial_u}\\
    \bar{u}_{k+j \mid k} &\in\bar{\mathbb{U}}_j,\ \ \ \ \; j=0,1,\dots,N-1, \label{subeq:input_r} \\
    \hat{u}_{k+j \mid k} &\in\bar{\mathbb{U}}_j,\ \ \ \ \; j=0,1,\dots,N-1, \label{subeq:input_l} \\
   \bar{x}_{k+j \mid k} & \in\bar{\mathbb{X}}_j, \ \ \ \ \; j=0,1,\dots,N, \label{subeq:state_r} \\
   \hat{x}_{k+j \mid k} & \in\hat{\mathbb{X}}_{j \mid k}, \ \ \ j=0,1,\dots,N. \label{subeq:state_l}
\end{align}
\end{subequations}
The OCP features two separate input sequences (also referred to as two \textit{horizons}), namely a robust input sequence \mbox{$\bar{U}_k = (\bar{u}_{k \mid k},\bar{u}_{k+1 \mid k},\dots,\bar{u}_{k+N-1 \mid k})$} used in a nominal prediction model \eqref{subeq:fr} of $\mathcal{S}$ with $w_k=0, k\in\mathbb{N}$, and the performance input sequence $\hat{U}_k = (\hat{u}_{k \mid k},\hat{u}_{k+1 \mid k},\dots,\hat{u}_{k+N-1 \mid k})$ used in the LB prediction model \eqref{subeq:fL} of $\mathcal{S}$. 
Here, the notation $u_{k + j \mid k}$ is used to denote the prediction of the variable $u$ at time $k+j$ made at time step $k$.
The LB prediction model in \eqref{subeq:fL} consists of the nominal dynamics $f$ and the residual dynamics estimated by $d(\hat{x}_{k+j \mid k},\hat{u}_{k+j \mid k},\mathcal{D}_k)$,
using a data set $\mathcal{D}_k$ available at time step $k$. The data set $\mathcal{D}_k$ can consist of past collections of state and input samples. Note that in various applications this data is not available a priori and, hence, the model must be adapted online. If no data are available a priori, then $\mathcal{D}_0 = \emptyset$.
The two horizons are coupled by sharing the initial state \eqref{subeq:same_initial_x}, i.e., $\bar{x}_{k \mid k} = \hat{x}_{k \mid k}$ and the first control input \eqref{subeq:same_initial_u}, i.e., $\bar{u}_{k \mid k} = \hat{u}_{k \mid k}$.
Note that in \eqref{subeq:same_initial_x} the initial state $\hat{x}_{k \mid k} = \bar{x}_{k \mid k}$ must lie in the Minkowski sum $x_k \oplus X = \{ x_k + \xi \mid \xi \in X \}$, where $X$ can be a robust positive invariant set for $\mathcal{S}$, as sometimes used in tube-based MPC \cite{bib:book_MPC}, or $X = \{0\}$ as used in classical robust MPC.
Both input sequences are constrained according to \eqref{subeq:input_r} and \eqref{subeq:input_l}, where the (possibly tightened) input constraint sets $\bar{\mathbb{U}}_j \subseteq \mathbb{U}$ can be selected as, e.g., in tube-MPC \cite{bib:book_MPC}.
Similarly, we consider (tightened) constraint sets $\bar{\mathbb{X}}_j\subseteq\mathbb{X},j=0,1,\dots,N$ in \eqref{subeq:state_r} as is common in various RMPC approaches \cite{bib:geurts_robust,bib:book_MPC}.
The LB horizon uses the data set $\mathcal{D}_k$ to adapt the tightened state constraint set $\hat{\mathbb{X}}_{j \mid k}\subseteq\mathbb{X} $ in \eqref{subeq:state_l} over time $k \in \mathbb{N}$.
Finally, let $(\bar{x}^*_{k \mid k},\bar{U}_k^*,\hat{U}_k^*)$ be a minimizer of \eqref{eq:MPC_problem} for $x_k$, which we assume exists, at time $k\in\mathbb{N}$ for cost function $J$.
For simplicity we assume that the minimizer is unique, but all the results hold also in case of non-uniqueness.
 
The CMPC can be perceived as a combination of two separate MPCs, namely, an RMPC
\begin{subequations}
\label{eq:RMPC_problem}
\vspace{-2pt}
\begin{align}
  \min_{\bar{x}_{k \mid k},\bar{U}_k}&~ J_R(x_k, \bar{x}_{k \mid k},\bar{U}_k)\\
  \text{s.t. }&~  \eqref{subeq:fr},\ \eqref{subeq:input_r},\ \eqref{subeq:state_r}, \ \bar{x}_{k \mid k} \in x_{k} \oplus X \label{subeq:rmpc_state}
\end{align}
and an LB-MPC
\end{subequations}
\begin{subequations}
\label{eq:LBMPC_problem}
\vspace{-4pt}
\begin{align}
  \min_{\hat{x}_{k \mid k},\hat{U}_k}&~ J_L(x_k, \hat{x}_{k \mid k},\hat{U}_k)\\
  \text{s.t. }&~  \eqref{subeq:fL},\ \eqref{subeq:input_l},\ \eqref{subeq:state_l}, \ \hat{x}_{k \mid k} \in x_{k} \oplus X
\end{align}
\end{subequations}
that are coupled by \eqref{subeq:same_initial_x} and \eqref{subeq:same_initial_u}. The cost function $J$ can be chosen, for instance, as a convex combination of $J_R$ and $J_L$.
\begin{figure}[t]
    \vspace{4pt}
    \centering\includegraphics[width=\columnwidth]{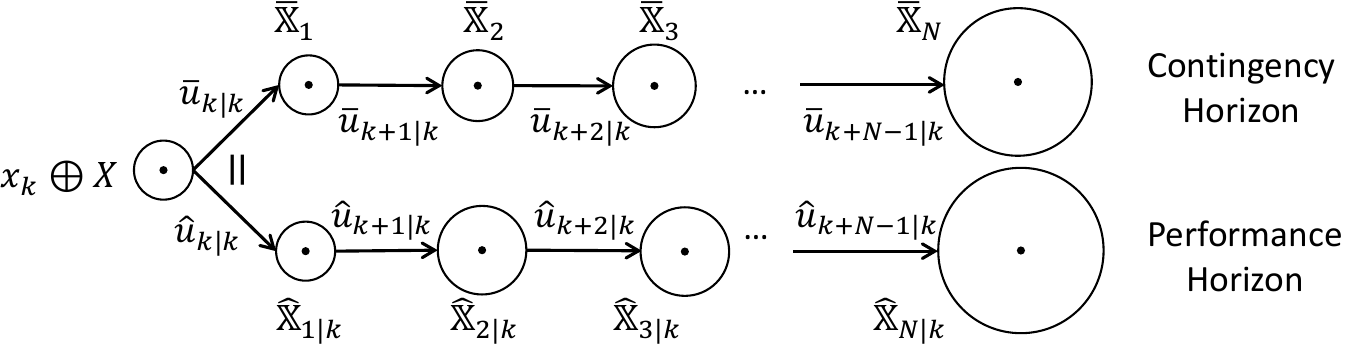}
    \vspace{-14pt}
    \caption{A schematic of CMPC framework with circles indicating the admissible state constraint sets for each horizon.}
    \label{fig:CMPC_scheme}
    \vspace{-14pt}
\end{figure}
The CMPC \eqref{eq:MPC_problem} can balance the performance of the LB-MPC, while ensuring safety through the RMPC under appropriate assumptions, as we will demonstrate in this paper.
Fig. \ref{fig:CMPC_scheme} depicts the CMPC with the robust \textit{contingency horizon} and the learning-based \textit{performance horizon}.
Note that CMPC \eqref{eq:MPC_problem} can accommodate various RMPC and LB-MPC methods by appropriately designing its components \eqref{eq:RMPC_problem} and \eqref{eq:LBMPC_problem}.

The resulting control input applied to $\mathcal{S}$ %
at time $k$ is $u_k = \pi(x_k) := \, \kappa(x_k, \bar{x}^*_{k \mid k}) +\bar{u}_{k \mid k}^{*} = \kappa(x_k, \hat{x}^*_{k \mid k}) + \hat{u}_{k \mid k}^{*}$, where an ancillary control law $\kappa$ is used as, e.g., in some tube-based MPC setups  \cite{bib:book_MPC}.
Note that for classical robust MPC, where $X = \{0\}$, the ancillary control law $\kappa$ is only used for constraint tightening and is not used in the closed-loop control law as $\kappa$ is typically zero for $x_k = \bar{x}_{k \mid k}= \hat{x}_{k \mid k}$.
In tube-based MPC the nominal state is relaxed to lie in a robust positive invariant set $X$ around the nominal trajectory \cite{bib:book_MPC}.

\vspace{-2pt}
\section{Safety Guarantees}\label{sec:safety}
In this section, we present the assumptions on the adopted RMPC \eqref{eq:RMPC_problem} and LB-MPC \eqref{eq:LBMPC_problem} that ensure robust recursive feasibility of the CMPC \eqref{eq:MPC_problem}, i.e., if the CMPC is feasible at time $k\in\mathbb{N}$ then it is also feasible at the next time step $k+1$. 

\begin{assumption}
\label{assum:terminal_set}
The RMPC \eqref{eq:RMPC_problem} is robustly recursively feasible in the sense that, if there is a pair $(\bar{x}_{k \mid k}, \bar{U}_k)$ that is feasible for \eqref{eq:RMPC_problem}, for state $x_k$ at time $k \in \mathbb{N}$, then there is, for all $w_k\in\mathbb{W}$, a pair $(\bar{x}_{k+1 \mid k+1}, \bar{U}_{k+1})$ that is feasible for \eqref{eq:RMPC_problem}, for the state $x_{k+1} = f(x_k, \kappa(x_k, \bar{x}_{k \mid k}) +\bar{u}_{k \mid k}) + w_k$ at time $k+1$.
\end{assumption}
    
\begin{assumption}
\label{assum:disturbance}
The learned residual dynamics satisfy $d(x_k,u_k,\mathcal{D}_k) \in \mathbb{W}$, and it holds that $\hat{\mathbb{X}}_{j \mid k} \supseteq \bar{\mathbb{X}}_j$, $j\in\{0,1,\dots,N\}$, for all $k\in\mathbb{N}$.
\end{assumption}

Note that Ass. 1 can be guaranteed for various robust MPC schemes, such as tube-based MPC, using proper constraint tightening and terminal ingredients \cite{bib:alvarado_limon, bib:book_MPC}. 
Based on these assumptions, we can show that
the CMPC \eqref{eq:MPC_problem} inherits the robust recursive feasibility of the RMPC \eqref{eq:RMPC_problem}:

\begin{theorem}
\label{thm:RF_CMPC}
Let Ass. \ref{assum:terminal_set} and \ref{assum:disturbance} hold. 
If OCP \eqref{eq:MPC_problem} has a feasible solution $(\bar{x}_{k \mid k},\bar{U}_{k}, \hat{U}_{k})$ for $x_k$ at time $k\in\mathbb{N}$, then \eqref{eq:MPC_problem} is feasible for the next state $x_{k+1}=f(x_k, \kappa(x_k, \bar{x}_{k \mid k}) +\bar{u}_{k \mid k})+w_k$, for any $ w_k\in\mathbb{W}$, at time $k+1$.
\end{theorem}

\begin{proof}
Given a feasible initial nominal state $\bar{x}_{k \mid k}$ and a pair of input sequences $(\bar{U}_k,\hat{U}_k)$, for state $x_k$ at time $k\in\mathbb{N}$, we have that the corresponding predicted state trajectories satisfy $(\bar{x}_{k \mid k},\bar{x}_{k+1 \mid k},\dots,\bar{x}_{k+N \mid k})\in  ( [ x_k \oplus X ] \cap \bar{\mathbb{X}}_0 ) \times\bar{\mathbb{X}}_1 \times\dots\times\bar{\mathbb{X}}_{N-1}\times\bar{\mathbb{X}}_N$ and $(\hat{x}_{k \mid k},\hat{x}_{k+1 \mid k},\dots,\hat{x}_{k+N \mid k})\in  ( [x_k \oplus X ] \cap \hat{\mathbb{X}}_0 ) \times\hat{\mathbb{X}}_{1 \mid k} \times\dots\times\hat{\mathbb{X}}_{N-1 \mid k}\times\hat{\mathbb{X}}_{N \mid k}$.
Due to Ass. \ref{assum:terminal_set}, this implies that there exist a nominal initial state $\bar{x}_{k+1 \mid k+1}$ and an input sequence $\Bar{U}_{k+1}$ that satisfy \eqref{subeq:rmpc_state}
at time $k+1$ for the state $x_{k+1} = f(x_k, \kappa(x_k, \bar{x}_{k \mid k}) +\bar{u}_{k \mid k}) + w_k$. Thus, for the corresponding state sequence, we have that $(\bar{x}_{k+1 \mid k+1}, \bar{x}_{k+2 \mid k+1}, \dots, \bar{x}_{k+N \mid k+1}, \bar{x}_{k+N+1 \mid k+1}) \in  ([ x_{k+1} \oplus X ] \cap \bar{\mathbb{X}}_0 )  \times \bar{\mathbb{X}}_1\times \dots \times \bar{\mathbb{X}}_{N-1} \times \bar{\mathbb{X}}_N$.
It remains to show that there exist a feasible initial nominal state and an input sequence for the performance horizon of problem \eqref{eq:MPC_problem}. 
Due to Ass. \ref{assum:disturbance}, we have that $d(x_k,u_k,\mathcal{D}_{k+1}) \in \mathbb{W}$ and $\hat{\mathbb{X}}_{j \mid k+1}  \supseteq \bar{\mathbb{X}}_j$.
This implies that we can select $\hat{x}_{k+1 \mid k+1} := \bar{x}_{k+1 \mid k+1}$ and $\hat{U}_{k+1} := \bar{U}_{k+1}$ that leads to a state sequence $(\hat{x}_{k+1 \mid k+1},\hat{x}_{k+2 \mid k+1},\dots,\hat{x}_{k+N \mid k+1},\hat{x}_{k+N+1 \mid k+1})\in  ( [ x_{k+1} \oplus X \cap \hat{\mathbb{X}}_{0 \mid k+1} ] ) \times\hat{\mathbb{X}}_{1 \mid k+1}\times\dots\times \hat{\mathbb{X}}_{N-1 \mid k+1}\times \hat{\mathbb{X}}_{N \mid k+1}$.
Hence, there exists a feasible solution $(\bar{x}_{k+1 \mid k+1},\bar{U}_{k+1}, \bar{U}_{k+1})$ for $x_{k+1}$.
\end{proof}
Theorem \ref{thm:RF_CMPC} shows that under Ass. \ref{assum:terminal_set} and \ref{assum:disturbance}, the CMPC inherits the robust recursive feasibility of the RMPC.
When the estimated residual dynamics $d$ fail to satisfy Ass. \ref{assum:disturbance}, we can use constraint softening for the LB-MPC to preserve robust recursive feasibility. Thereto, we replace \eqref{eq:MPC_problem} by
\begin{subequations}
\label{eq:MPC_problem_soft}
\begin{align}
  \min_{\bar{x}_{k \mid k}, \bar{U}_k,\hat{U}_k,\mathcal{E}_k}&~ J(x_k, \bar{x}_{k \mid k},\bar{U}_k,\hat{U}_k)+J_\varepsilon(\mathcal{E}_k)\label{subeq:soft_cost}\\
  \text{s.t. } \eqref{subeq:fr} & -\eqref{subeq:state_r}, \\
    \hspace{0mm} \hat{x}_{k+j \mid k} \in \hat{\mathbb{X}}_{j \mid k}&(\varepsilon_{k+j \mid k}), \ \ \varepsilon_{k+j \mid k}\geqslant0, \ j=0,1,\dots,N,\label{subeq:state_l_soft}
\end{align}
\end{subequations}
where $\mathcal{E}_k=(\varepsilon_{k \mid k},\dots,\varepsilon_{k+N \mid k}) \in \mathbb{R}^{N+1}_{\geq0}$ is a sequence of nonnegative slack variables and a corresponding penalty term $J_{\varepsilon} : \mathbb{R}^{N+1}_{\geqslant 0} \rightarrow \mathbb{R}_{\geqslant 0}$ is added to the cost function, with $J_\varepsilon(\mathcal{E})>0$, if $\mathcal{E}\neq0$, and $J_\varepsilon(0)=0$. 
The state constraints $\hat{\mathbb{X}}_{j \mid k}$ in \eqref{eq:MPC_problem} are assumed to be softened such that $\hat{\mathbb{X}}_{j \mid k}(0)=\hat{\mathbb{X}}_{j \mid k}$ and  $\hat{\mathbb{X}}_{j \mid k}(\varepsilon)\rightarrow\mathbb{R}^{n_x}$ when $\varepsilon \rightarrow\infty$.
\begin{theorem}
\label{th:rf_soft}
Let Ass. \ref{assum:terminal_set} hold. If $(\bar{x}_{k \mid k},\bar{U}_{k}, \hat{U}_{k}, \mathcal{E}_k)$ is a feasible quadruple for OCP \eqref{eq:MPC_problem_soft} for $x_k$ at time $k\in\mathbb{N}$, then \eqref{eq:MPC_problem_soft} is feasible for the state $x_{k+1}=f(x_k,\kappa(x_k, \bar{x}_{k \mid k}) +\bar{u}_{k \mid k})+w_k$ for any $ w_k\in\mathbb{W}$ at time $k+1$.
\end{theorem}

\begin{proof}
    Similar to the proof of Thm. \ref{thm:RF_CMPC}, based on Ass. \ref{assum:terminal_set}, there exist a feasible nominal initial state $\bar{x}_{k+1 \mid k+1}$ and an input sequence $\Bar{U}_{k+1}$ that satisfies \eqref{subeq:rmpc_state} at time $k+1$ for $x_{k+1} = f(x_k, \kappa(x_k, \bar{x}_{k \mid k}) +\bar{u}_{k \mid k}) + w_k$.
    Due to $\lim_{\varepsilon \rightarrow \infty} \hat{\mathbb{X}}_{j \mid k}(\varepsilon) = \mathbb{R}^{n_x}$ for $k \in \mathbb{N}$, there always exists a feasible $\mathcal{E}_{k+1} \geqslant 0$ large enough such that $\bar{\mathbb{X}}_j \subseteq \hat{\mathbb{X}}_{j \mid k+1}(\varepsilon_{k+1+j \mid k+1})$.
    It follows directly from the proof of Thm. 1 that $\hat{x}_{k+1 \mid k+1} := \bar{x}_{k+1 \mid k+1}$ and $\hat{U}_{k+1} := \bar{U}_{k+1}$ satisfy \eqref{eq:MPC_problem_soft}.
    Hence, there exists a feasible solution $(\bar{x}_{k+1 \mid k+1},\bar{U}_{k+1}, \bar{U}_{k+1}, \mathcal{E}_{k+1})$ for $x_{k+1}$.
\end{proof}

Note that the theoretical results of this general CMPC framework can be leveraged for various choices of the RMPC and LB-MPC as \eqref{eq:MPC_problem} -- \eqref{eq:MPC_problem_soft} are formulated on a high, conceptual level. Below, we demonstrate a specific implementation of CMPC for safe learning. Firstly, we introduce a relevant case study in the field of automated driving.

\section{Automated Lane Merging Application}
\label{sec:application}
\begin{figure}[t]
    \centering
     \includegraphics[trim={0 0.2cm 0 0},clip,width=0.8\columnwidth]{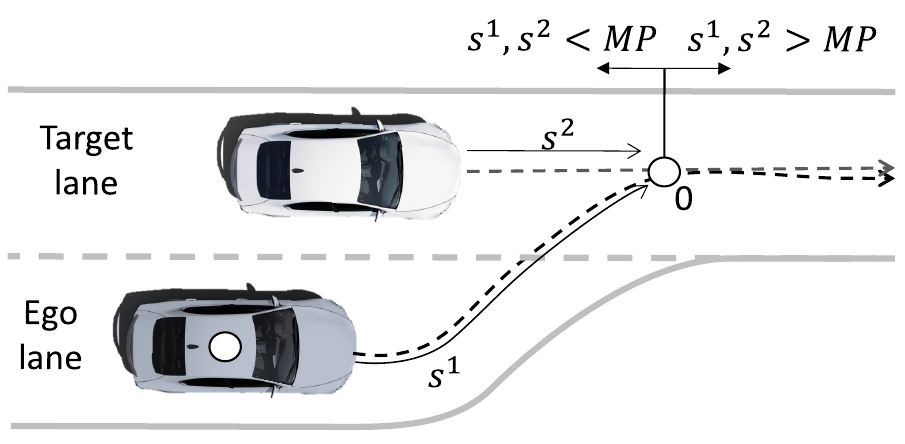}
    \caption{Set-up of the lane-merging scenario.}
    \label{fig:scenario}
    \vspace{-12pt}
\end{figure}

We consider a lane merging application, in which the ego vehicle should execute the following scenario, depicted in \mbox{Fig. \ref{fig:scenario}}:
\mbox{(i) Initially, the} ego vehicle (Agent 1, gray vehicle) is driving in the ego lane. The target vehicle (Agent 2, white vehicle) is driving in the target lane, driving with an unknown but bounded acceleration.
(ii) The ego lane is closing off, requiring \mbox{Agent 1} to merge into the target lane, before the merging point \textit{MP}.
(iii) After the merge, Agent 1 maintains a safe distance to Agent 2. When Agent 1 is in front, Agent 2 maintains a safe distance to Agent 1.
During the scenario, it is assumed that the paths of \mbox{Agent 1} and 2 are fixed (but not the speed along the path);
and the current position and velocity of \mbox{Agent 2} are known,  but its current and future acceleration profile are unknown.

The dynamics of \mbox{Agent 1} and 2 are formulated in terms of a one-dimensional path coordinate $s^i_k$ and the forward velocity $v^i_k$ for $i\in\{1,2\}$, as well as the relative distance \mbox{$\Delta s_k = s^2_k-s^1_k$} and the relative velocity $\Delta v_k = v^2_k-v^1_k$ of Agents 1 and 2.
Here, $s^i_k < 0$ and $s^i_k > 0$ indicate a position before and after the \textit{MP}, respectively, as depicted in Fig. \ref{fig:scenario}.
The acceleration of \mbox{Agent 1} $u_k^1$ is controlled by the CMPC, while the acceleration of \mbox{Agent 2} $u_k^2$ is modeled as an unknown, bounded disturbance.
The discrete time $k\in\mathbb{N}$ corresponds to real time $t=kT_s$ with the sampling period $T_s$. 
We select the state as $x_k = \begin{bmatrix} \Delta s_k~\Delta v_k~ s_k^1~v_k^1 \end{bmatrix}^\top\in\mathbb{R}^4$ with
\begin{equation}\label{eq:dynamics_discrete}
    x_{k+1}^{} =A^{}x_{k}^{}+B_1u_{k}^{1}+B_2u_k^2,\quad\text{where}
\end{equation}
\begin{equation*}
{\small A = \begin{bmatrix} 1&T_s&0&0\\0&1&0&0\\0&0&1&T_s\\0&0&0&1 \end{bmatrix} \text{, } B_1^{} = \begin{bmatrix}-\frac{1}{2}T_s^2\\-T_s\\\frac{1}{2}T_s^2\\T_s\end{bmatrix} \text{\normalsize and } B_2^{} = \begin{bmatrix}\frac{1}{2}T_s^2\\T_s\\0\\0\end{bmatrix}.} 
\end{equation*}
 We assume that the process noise $v_k = 0$ and the bounded disturbance $w_k=B_2u^2_k = g(x_k) \in \mathbb{W}$, $u^2_k \in [ u_{\text{min}}^2, u_{\text{max}}^2 ]$ is a state-dependent policy for Agent 2 and $ v_{\text{min}}^2\leqslant v_k^2\leqslant v_{\text{max}}$ for $k \in \mathbb{N}$.
 The acceleration and velocity of \mbox{Agent 1} are bounded, yielding the input constraint set
$\mathbb{U}:= [ u_{\text{min}}^1,  u_{\text{max}}^1]$ and
$0\leqslant v^1_k\leqslant  v_{\text{max}}$ for all $k\in\mathbb N$.
Finally, there is a smooth state-dependent safety distance function $D_{\text{safe}}$, which depends on the relative distance $ \Delta s_k$, velocity $ v_k^1$ and position $ s_k^1$ of \mbox{Agent 1}, that ensures a safe distance to Agent 2:
\begin{equation}\label{eq:safety_constraint}
     D_{\text{safe}}( s^1_k,\Delta s_k,v^1_k ) \leqslant 0, \quad k\in\mathbb{N}.
\end{equation}
We refer to \cite{bib:geurts_robust} for details on the design of the safety distance function  $D_{\text{safe}}$. 
This leads to the state constraint set
 \begin{equation}
\label{eq:state_constraint}
    \mathbb{X}:=\{ x\in\mathbb{R}^4 \mid 0\leqslant v^1\leqslant v_{\text{max}}, \, D_{\text{safe}}(s^1,\Delta s,v^1) \leqslant 0 \}.
\end{equation}
The objective is to use CMPC to control Agent 1, completing steps (i) to (iii) of the scenario, learning the acceleration policy $g$ of Agent 2 and adhering to all the safety constraints.

\section{CMPC for Automated lane merging}
\label{sec:implementation}
 In this section, we apply the CMPC framework proposed in Sec. \ref{sec:framework} to the automated lane merging problem (see Sec. \ref{sec:application}) by combining the RMPC from \cite{bib:geurts_robust} for safety guarantees and a GP-MPC as in \cite{bib:baltussen} for performance.
We use the soft-constrained LB-MPC variant of CMPC described in \eqref{eq:MPC_problem_soft}.

\subsection{Robust MPC - Contingency Horizon}
 For the contingency horizon, we use a classical robust MPC with  a nominal prediction model \eqref{subeq:fr} as
\begin{equation}\label{eq:dynamics_robust_mpc}
    \bar{x}_{k+1 \mid k}^{} = f(\bar{x}_{k \mid k}^{},\bar{u}_{k \mid k}^{1}) = A^{}\bar{x}_{k \mid k}^{}+B_1\bar{u}_{k \mid k}^{1},
\end{equation}
predicted state $\bar{x}_{k \mid k} := [ \Delta \bar{s}_{k\mid k} ~ \Delta \bar{v}_{k\mid k} ~ \bar{s}_{k\mid k}^1 ~ \bar{v}_{k\mid k}^1 ]^\top$
  and $\kappa(x,\bar{x})=0$ for all $x$, $\bar{x}$ and $X = \{0\}$.
We construct the tightened constraint set $\bar{\mathbb{X}}_j := \mathbb{X}\ominus \sum_{i=0}^{j-1}A^i \mathbb{W}$,  for $j=1,2,\dots,N-1$ in \eqref{subeq:state_r} (see details in \cite{bib:geurts_robust}),
 where the Pontryagin set difference of $\mathbb{A}\subseteq \mathbb{R}^n$ and $\mathbb{B}\subseteq\mathbb{R}^n$ is $\mathbb{A}\ominus\mathbb{B}:=\{a \mid \forall b\in \mathbb{B}, a+b\in \mathbb{A} \}$.
 We consider two disjoint robust control invariant sets $\Omega_1$ (for merging behind Agent 2) and $\Omega_2$ (for merging in front of Agent 2) from \cite{bib:geurts_robust}.
The terminal constraint set of the contingency horizon $\bar{\mathbb{X}}_N$ in \eqref{subeq:state_r} is defined as the union of $\Omega_1$ and $\Omega_2$, such that $\bar{\mathbb{X}}_N = \Omega_1 \cup \Omega_2$.
In \cite{bib:geurts_robust}, the details of $\Omega_1$ and $\Omega_2$ are given 
and it is shown that the adopted constraint tightening and choice for $\bar{\mathbb{X}}_N$ leads to the robust recursive feasibility of the RMPC as in Ass. \ref{assum:terminal_set}.

\subsection{Gaussian Process MPC - Performance Horizon}
The prediction model \eqref{subeq:fL} of the  performance horizon consists of the nominal dynamics $f$ as in \eqref{eq:dynamics_robust_mpc} and augmented dynamics $d$ to account for the unmodeled dynamics $w_k = B_2u_k^2 = g(x_k)$,  using GP regression \cite{bib:rasmussen_gp} and Bayesian inference based on collected observations of the uncertainty
\begin{equation}
    y_k = B_2^{\dagger} \bigl( x_{k+1} - f(x_{k},u_k^1) \bigr) = u^2_k,
\end{equation}
and the regressor $\mathbf{z}_k \hspace{-0mm} = \hspace{-0mm} x_k$,  where $B_2^\dagger$ denotes the Moore-Penrose pseudo-inverse of $B_2$.
We use a data set $\mathcal{D} = \left\{\mathbf{Z}, \mathbf{y} \right\}$ with the collected training data 
$\mathbf{Z} = [\mathbf{z}_1, \mathbf{z}_2, \dots, \mathbf{z}_{n_D}] \in \mathbb{R}^{n_z \times n_D}$ and ${\mathbf{y} = [y_1, y_2, \dots, y_{n_D} ] \in \mathbb{R}^{1\times n_D}}$.
 At each time $k$, we append the data set with the collected training data: $\mathcal{D}_{k+1} = \mathcal{D}_k \cup \{\mathbf{z}_{k}, \mathbf{y}_k\}$.
We use a GP $\Tilde{d} : \mathbb{R}^{n_z} \rightarrow \mathbb{R}$ to approximate the policy of Agent 2:
\begin{equation}\label{eq:dynamics_GP_mpc}
    B_2^{\dagger} g(x) \approx \Tilde{d}\left(\mathbf{z} \right) \sim \mathcal{N} \bigl(d\left(\mathbf{z} \right), \Sigma^{d}\left(\mathbf{z} \right) \bigr),
\end{equation}
where $d(\mathbf{z})$ and $\Sigma^{d}(\mathbf{z})$ denote the mean and covariance function of the posterior of the GP at a given test point $\mathbf{z}$:
\begin{equation}
    d(\mathbf{z}) = \mathbf{k}_{\mathbf{z Z}} K^{-1}_{\mathbf{ZZ}} \mathbf{y},~ 
    \Sigma^d(\mathbf{z}) =  k\left(\mathbf{z}, \mathbf{z}\right) - \mathbf{k}_{\mathbf{z Z}} K^{-1}_{\mathbf{Z}\mathbf{Z}} \mathbf{k}_{\mathbf{Z}\mathbf{z}}
\label{eq:Sparse_GP}
\end{equation}
where $k$ denotes the user-defined kernel function, {$\mathbf{k}_{\mathbf{Z z}} \in \mathbb{R}^{n_D}$} is the concatenation of the kernel function evaluated at the test point {$\mathbf{z}$} and the training set {$\mathbf{Z}$}, where {$[\mathbf{k}_{\mathbf{Z z}}]_i = k\left(\mathbf{z}_i,\mathbf{z}\right)$} and {${\mathbf{k}_{\mathbf{Z z}}^\top = \mathbf{k}_{\mathbf{z Z}}}$}, and $K_{\mathbf{ZZ}} \in \mathbb{R}^{{n_D} \times {n_D}}$ is a Gram matrix which satisfies {${[K_{\mathbf{Z Z'}}]_{ij} = k\bigl(\mathbf{z}_i,\mathbf{z}'_{j}\bigr)}$}.
In this work, we employ a squared exponential kernel function:
\begin{equation}
\label{eq:SE_ker}
    \hspace{-3mm} k\left(\mathbf{z}, \mathbf{z}'\right) = \sigma^2_d \exp \bigl( - \tfrac{1}{2} \left(\mathbf{z} \hspace{-0.3mm} -\hspace{-0.3mm} \mathbf{z}' \right)^\top \hspace{-0.4mm} L_d^{-2} \left(\mathbf{z} \hspace{-0.3mm}-\hspace{-0.3mm} \mathbf{z}'\right) \bigr), \hspace{-2mm}
\end{equation}
 where $\sigma_d$ denotes the prior covariance and $L_d$ is the length-scale matrix.
The GP formulation in \eqref{eq:dynamics_GP_mpc} is generally intractable for multi-step predictions. To overcome this issue, we make several commonly used approximations in GP-MPC.  Similarly to \cite{bib:hewing}, (i) we assume that the predicted state $\hat{x}$ and the posterior GP $d$ are jointly Gaussian, and \mbox{(ii) the} consecutive GP evaluations are independent. In addition, (iii) we approximate the consecutive means, covariances and cross-covariances of the GP with a first-order Taylor approximation. 
 (iv) We use a Sparse Pseudo-Input GP with $M$ inducing points that are equally spaced over the prediction horizon. For further details, we refer the reader to \cite{bib:baltussen}.
The resulting prediction model \eqref{subeq:fL} of the  performance horizon is composed of the mean and covariance of the posterior GP: 
\begin{subequations}
\begin{align}
\hspace{-2mm} \hat{x}_{k+j+1 \mid k} & = A \hat{x}_{k+j \mid k}  + B_1 \hat{u}^1_k+ B_2 d\bigl( \hat{x}_{k+j \mid k}, \mathcal{D}_k \bigr), \label{eq:Model_Mean} \\
\hspace{-2mm} \Sigma_{k+j+1 \mid k}^{x} & =  \begin{bmatrix}
     A &  B_2
\end{bmatrix} \Sigma_{k+j \mid k}  \begin{bmatrix}
    A &  B_2
\end{bmatrix}^{\top},
\label{eq:Model_Cov}
\end{align}
\label{eq:Sparse_GP_model}
\end{subequations}
for {$j = 0, 1,\dots,N-1$},
where the initial prediction equals the current state {$\hat{x}_{0 \mid k} = x_k$} with {${\Sigma^{x}_{ 0 \mid k} = \mathbf{0}}$}, and $\Sigma$ is the predicted joint covariance matrix of the state $x$ and  $\tilde{d}(\mathbf{z})$ \cite{bib:baltussen}.

We use the GP prediction model to obtain an estimate of the lower and upper bound of $\Delta s$ based on the predicted sequence $\Hat{U}_k$ and the observed data set $\mathcal{D}_k$.
In the CMPC \eqref{eq:MPC_problem_soft}, we predict the future states through the mean function \eqref{eq:Model_Mean} and use the state covariance \eqref{eq:Model_Cov} to adapt  the set $\hat{\mathbb{X}}_{j \mid k}$.
This constraint is softened in the performance horizon by an $\ell_1$ penalty method:
\begin{align}
    &\hat{\mathbb{X}}_{j \mid k}:=\{ x_{k+j  \mid  k} \in\mathbb{R}^4 \mid 0\leqslant v^1_{k+j  \mid  k} \leqslant v_{\text{max}}, \\ 
    &D_{\text{safe}}(s^1_{k+j  \mid  k},\Delta s_{k+j  \mid  k} \pm 2\sigma_{s^2,k+j  \mid  k},v^1_{k+j  \mid  k}) \leqslant \varepsilon_{k+j \mid k} \} \nonumber,
 \end{align}
where $\sigma_{s^2,k+j}$ denotes the approximate standard deviation of the position of Agent 2 at prediction step $k+j$.

\subsection{Integration into CMPC Framework}
We choose the  primary cost function  $J$ in \eqref{subeq:soft_cost} as
\begin{align}
       J&(x_k,  \bar{x}_{k \mid k},\bar{U}_k,\hat{U}_k) := P H(x_k,\bar{U}_k)+(1-P) H(x_k,\hat{U}_k), \nonumber \\
        H&({x}_k,{U}_k) = \hspace{-1mm} \sum_{j=0}^{N}
        \|{ v_{\text{ref}}^1}-{v}^1_{k+j \mid k}\|_Q^2 + \hspace{-1mm}\sum_{j=0}^{N-1} \|{u}^1_{k+j \mid k}\|_R^2 \\ 
        & + \|\Delta u^1_{k+j \mid k}\|_S^2, \nonumber
\end{align}
where $\|z\|_Q^2 = z^\top Q z$, and $Q,R,S>0$ are tuning parameters, $\Delta u^1_{k+j \mid k} := u^1_{k+j \mid k}-u^1_{k+j-1 \mid k}$ is the change in control input with \hspace{0.5mm}$u^1_{k-1 \mid k} := u^1_{k-1}$, and $v_{\text{ref}}^1$ denotes the desired velocity of Agent 1.
The parameter $P\in[0,1]$ balances the cost for the contingency and the performance horizon.
Further, we use an $\ell_1$ penalty function $J_\epsilon (\mathcal{E}) = \rho \|\mathcal{E}\|_1$ in \eqref{subeq:soft_cost} with a penalty weight $\rho > 0$.
While the GP-MPC is soft constrained over the performance horizon, the recursive feasibility and closed-loop safety of the CMPC are ensured by Thm. \ref{th:rf_soft}.

\section{Numerical Results}\label{sec:experiments}
\begin{figure*}[h]
    \centering
    \includegraphics[width=1\textwidth]{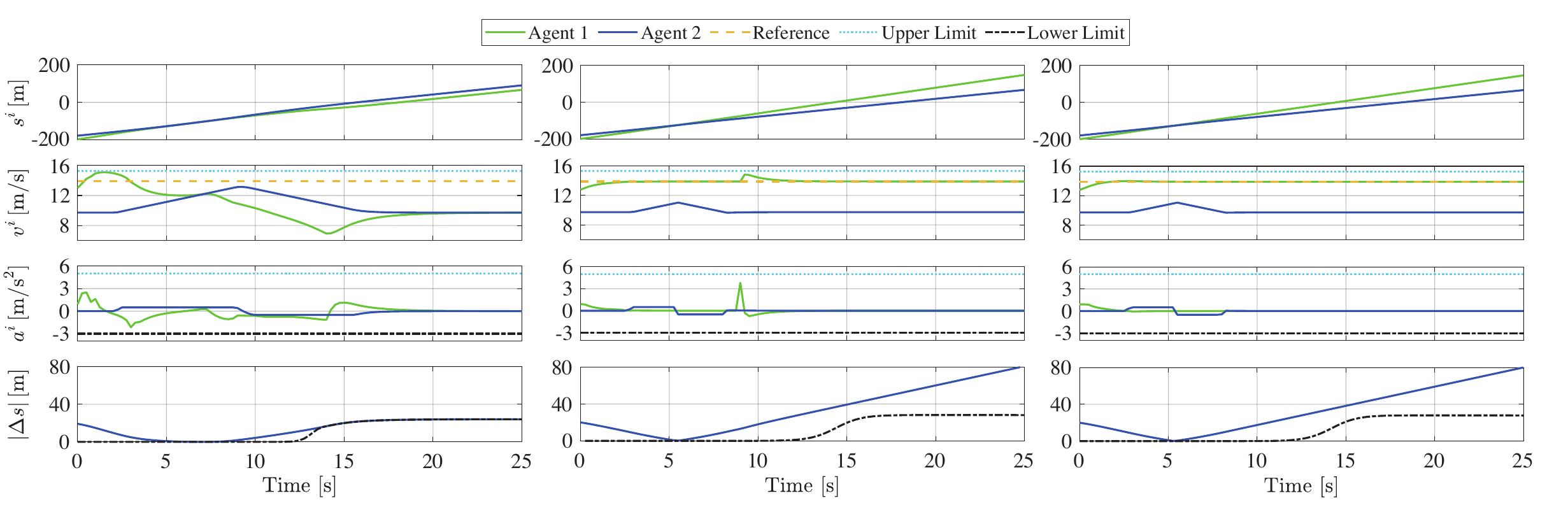}
    \vspace{-16pt}
    \caption{Comparison of RMPC (left), CMPC (middle) and GP-MPC (right) for one lane merging scenario, see Sec. \ref{sec:use_case}. Shown are the trajectories of the path coordinate $s^i$, velocity $v^i$, acceleration $u^i$ and relative distance $|\Delta s|$ of \mbox{Agent 1} and 2. The lower limit of $|\Delta s|$ is imposed by $D_{\text{safe}}$.}
    \label{fig:CMPC_merge_comparison}
    \vspace{-9pt}
\end{figure*}

The behavior of the CMPC detailed in Sec. \ref{sec:application}-\ref{sec:implementation} is verified in a numerical case study, which compares CMPC with the lane-merging controllers using standalone RMPC and GP-MPC, and shows the benefits of CMPC.
The following parameters are used in the case study\footnote{The numerical experiments are performed in Matlab R2023a with IPOPT 3.14.11 \cite{bib:ipopt} and CasADi 3.6.5 \cite{bib:casadi} on an Intel Core i7-9750H at 2.60 $\mathrm{GHz}$.\label{fn:computation}}: $v_{\text{ref}}^1 = \frac{50}{3.6}\mathrm{\ \frac{m}{s}}$, $v_{\text{max}} = 1.1\cdot v_{\text{ref}}^1$, $T_s = 0.25 \mathrm{\ s}$, $N = 20$, 
$Q=10$, $R=1$, $S = 10$, $\rho = 10^4$ $u_{\text{min}}^1 = -3 \mathrm{\ \frac{m}{s^2}}$, $u_{\text{max}}^1 = 5\mathrm{\ \frac{m}{s^2}}$, $u_{\text{min}}^2=-0.5\mathrm{\ \frac{m}{s^2}}$, $u_{\text{max}}^2=0.5\mathrm{\ \frac{m}{s^2}}$, $v_{\text{min}}^2 = \frac{25}{3.6} \mathrm{\ \frac{m}{s}}$, $s_0^{1} = -200 \mathrm{\ m}$, $s_0^{2} = -180 \mathrm{\ m}$, $P=0.5$, $M=4$, $\sigma_d = 0.7$ and $L_d = \mathrm{diag}(5, 100, 500, 100)$. Parameters related to the safety function $D_{\text{safe}}$ are the same as in \cite{bib:geurts_robust}.
The GP-MPC starts with $n_D = 0$ data points and the size of $\mathcal{D}_k$ increments at every time step.
\mbox{Agent 2}'s acceleration $u_k^2$ adheres to the state-dependent cooperative driving policy
\begin{equation}\label{eq:cooperative_u2}
    u_k^2 = \min(\max(k_v(v_{\text{ref}}^2-v^2_k)+u_{\Delta s},u_{\text{min}}^2),u_{\text{max}}^2),
    \vspace{-2pt}
\end{equation}
with a deceleration dependent on the relative distance where
\begin{align*}
    & u_{\Delta s} = \begin{cases}
    u_{\Delta s}^*,&\text{if}\quad s^2_k\geqslant -200,~ s^1_k<0,\\
    0,&\text{otherwise},
    \end{cases} \\
   & u_{\Delta s}^* = \begin{cases}
        k_{\Delta s}(\Delta s_{\text{ref}}-\Delta s_k),&\text{if}\quad\Delta s_{\text{react}} \geqslant \Delta s_k\geqslant0,\\
        k_{\Delta s}(-\Delta s_{\text{ref}}-\Delta s_k),&\text{if}\quad-\Delta s_{\text{react}} \leqslant \Delta s_k\leqslant0,\\
        0,&\text{otherwise},
    \end{cases}
    \vspace{-2pt}
\end{align*}
with a reference velocity $v_{\text{ref}}^2 \hspace{-0.1mm} = \hspace{-0.1mm} v_0^2$, reference following distance $\Delta s_{\text{ref}} = 10 \ \mathrm{m}$, reaction distance $\Delta s_{\text{react}} = 10 \ \mathrm{m}$, velocity gain $k_v = 1.0461$ and a relative distance gain $k_{\Delta s} = 0.4472$.
\vspace{-7pt}

\subsection{Merging scenario comparison}
\label{sec:use_case}
We consider a scenario in which \mbox{Agent 1} is driving at $v_0^{1}= \frac{46}{3.6}\mathrm{\frac{m}{s}}$ behind \mbox{Agent 2}, driving at $v_0^{2}=\frac{35}{3.6}\mathrm{\frac{m}{s}}$.
The results are shown in Fig. \ref{fig:CMPC_merge_comparison}.
In all MPCs, \mbox{Agent 1} first accelerates towards the desired velocity $v_{\text{ref}}^1$.
Then, \mbox{Agent 1} needs to decide to merge in front or behind \mbox{Agent 2}.
The LB-MPCs use an advanced prediction model to have a better understanding of the environment compared to the nominal model of the RMPC.
The proposed CMPC framework leverages this LB-MPC to make more informed decisions through the cost function (i.e., select better control inputs) than the standalone RMPC.
In case of GP-MPC and CMPC, the GP estimates based on past observations that \mbox{Agent 2} will decelerate when \mbox{Agent 1} tries to merge in front.
Consequently, the GP-MPC and CMPC policies control \mbox{Agent 1} to merge in front.
Hence, the GP-MPC and CMPC show more assertive behavior and improved performance over the standalone RMPC.
Around $9.25$ s, the CMPC slightly accelerates as the \textit{MP} enters the prediction horizon, to ensure a sufficient safety distance to Agent 2 and retain the feasibility of the CMPC, see Fig. \ref{fig:CMPC_merge_comparison}.
However, the RMPC does not consider the interactions with Agent 2 and shows more conservative behavior, resulting in a merge behind Agent 2.
All MPCs maintain sufficient distance towards \mbox{Agent 2}.
Finally, \mbox{Agent 1} converges to either $v_{\text{ref}}^1$ or to \mbox{Agent 2}'s velocity.
In summary, the CMPC is less conservative than RMPC by leveraging the GP while, as opposed to GP-MPC, preserving robust recursive feasibility and thus robust safety guarantees.
The CMPC results are obtained with a maximum computation time of $2.517 \ \mathrm{s}$ and average computation time of $0.699 \ \mathrm{s}$ per MPC iteration on a standard laptop PC\footref{fn:computation}, showing that real-time tractability is reachable with code optimization and high-performance hardware.
 \vspace{-7pt}
\subsection{Sensitivity analysis}

We further verify the CMPC \eqref{eq:MPC_problem_soft} by comparing it to the standalone RMPC and GP-MPC with respect to safety and performance in a larger sensitivity analysis, based on the following key performance indicators (KPIs) averaged over $n = 231$ scenarios and $m = 161$ time steps:
\begin{enumerate}
    \item \textbf{Slack:} The amount of constraint violation measured by the slack over the complete performance horizon $\frac{1}{n}\sum_{i=1}^n \frac{1}{m}\sum_{k=0}^{m-1} \|\mathcal{E}_{k}^i \|_1$.
    \item \textbf{Cost:} The merging progress and effort measured by the closed-loop stage cost $\frac{1}{n} \sum_{i=1}^{n} \frac{1}{m} \sum_{k=0}^{m-1} Q{(v_{\text{ref}}^1}-v_{k,i}^1)^2
    +R(u_{k,i}^1)^2 + S (\Delta u^1_{k,i})^2$, with $u^1_{-1,i} := 0$.
    \item \textbf{Merging time:} Performance measured with the merging time, the time when \mbox{Agent 1}'s position is $s^1>0 \ \mathrm{m}$.
    \item \textbf{Result:} Number of scenarios in which \mbox{Agent 1} merges in front or behind \mbox{Agent 2} as a measure of assertiveness.
\end{enumerate}

In the sensitivity analysis, the CMPC, RMPC and GP-MPC are tested for the parameterization defined at the beginning of this section, \mbox{Agent 2}'s control policy in \eqref{eq:cooperative_u2} and initial conditions $v^1_0 \in \{40,41,\dots,50 \}\mathrm{\frac{km}{h}}$ and $v^2_0 \in \{30,31,\dots,50\} \mathrm{\frac{km}{h}}$. 
The resulting KPIs are shown in Table \ref{tab:sensitivity_cooperative}.
In this study, the CMPC is more assertive than RMPC with 88 vs. 72 front merges (22\% increase) while having a reduced cost of 2.14 vs. 2.89 (26\% reduction) with comparable merging times.
The performance when merging behind Agent 2 is comparable for RMPC and CMPC.
However, in some cases, \mbox{Agent 1} optimistically tries to merge in front, but eventually has to merge behind to guarantee safety, causing the CMPC to incur a higher cost.
Both the RMPC and CMPC never exceed the safety constraint in closed loop due to their recursive feasibility properties.
The GP-MPC is the least conservative in the sense that it merges in front more often than RMPC and CMPC.
Similarly, when GP-MPC merges behind, the cost is lowest as it takes less caution than RMPC and CMPC.
However, the lack of safety guarantees leads to an increased use of soft constraints, and hence, constraint violation. The GP-MPC exceeds the safety constraint by at most $0.98$ m. However, no collision occurred due to the presence of a safety margin in $D_{\text{safe}}$ (see \cite{bib:geurts_robust}).
The CMPC combines the strengths of GP-MPC and RMPC, by learning the residual dynamics in \eqref{eq:system_S}, hence, improving performance while remaining safe.

Note that for 12 and 11 initial conditions for the RMPC and CMPC, respectively, the robust controlled invariant set $\bar{\mathbb{X}}_N$ is not reachable from the initial condition, and, hence, safety cannot be guaranteed under all possible disturbances.
As a consequence, these numerical experiments were infeasible and are therefore excluded from Table \ref{tab:sensitivity_cooperative}.
The soft-constrained GP-MPC can handle these scenarios by using the slack variables at the risk of a potential collision.
Providing formal safety guarantees may come at a cost of excluding certain initial conditions in which safety cannot be guaranteed. However, the CMPC framework can be used to analyze and flag these scenarios, either a priori or on-line.

\begin{table}[t]
\caption{KPIs for the results front (F) and behind (B). Note that RMPC and CMPC have 12 and 11 infeasible initial states, respectively, that are excluded of these KPIs.}
\vspace{-2pt}
 \label{tab:sensitivity_cooperative}
 \resizebox{1\columnwidth}{!}{%
 \begin{tabular}{|l||ll|ll|ll|ll|}
 \hline
 \multicolumn{1}{|l||}{\textbf{\shortstack{MPC \\ \,} }} &
   \multicolumn{2}{l|}{\textbf{\shortstack{Number [-] \\ \, }}} & 
   \multicolumn{2}{l|}{\textbf{\shortstack{Slack [-] \\ \, }}} &
   \multicolumn{2}{l|}{\textbf{\shortstack{Cost [-] \\ \, }}} &
   \multicolumn{2}{l|}{\textbf{\shortstack{\strut  Merging \\ Time {[}s{]} }}} 
   \\ \cline{1-9}
  \textbf{Result} & \multicolumn{1}{l|}{  F} &   B &
   \multicolumn{1}{l|}{  F} &
    B &
   \multicolumn{1}{l|}{  F} &
    B &
   \multicolumn{1}{l|}{   F} &
    B
   \\ \hhline{|=||=|=|=|=|=|=|=|=|}
    RMPC  & \multicolumn{1}{l|}{72} & 147 &
   \multicolumn{1}{l|}{0} & 0 & 
   \multicolumn{1}{l|}{2.89} & 62.61 & 
   \multicolumn{1}{l|}{14.42} & 16.88 \\ \hline
    CMPC  & \multicolumn{1}{l|}{88} & 132 &
   \multicolumn{1}{l|}{$10^{-7}$} & $10^{-6}$ & 
   \multicolumn{1}{l|}{2.14} & 65.72 & 
   \multicolumn{1}{l|}{14.82} & 16.94 \\ \hline
GP-MPC & \multicolumn{1}{l|}{153} & 78 &
   \multicolumn{1}{l|}{1.36} & $10^{-7}$ & 
   \multicolumn{1}{l|}{2.41} & 38.71 & 
   \multicolumn{1}{l|}{14.74} & 16.93 \\ \hline
 \end{tabular}%
 }
 \vspace{-10pt}
\end{table}

\section{Conclusions}\label{sec:conclusion}
To address the challenge of controlling safety-critical systems, this research focuses on achieving both assertive and flexible learning behavior while ensuring safety guarantees -- a concept often referred to as “safe learning.” To this end, it introduces a framework that combines general robust MPC (RMPC) algorithms with general learning-based MPC (LB-MPC) algorithms.
It is proven that this CMPC framework inherits the robust recursive feasibility properties of the RMPC while realizing improved performance by learning from data.
The proposed CMPC framework is applied to an automated driving use case of lane merging.
This case study shows that CMPC has improved driving performance due to the presence of the LB-MPC (i.e., merging in front more often, by learning the behavior of the other vehicle) compared to RMPC, while preserving the safety of the RMPC.

\end{document}